\documentclass{article}
\usepackage{amsfonts}

\usepackage{amsmath}

\newtheorem{theorem}{Theorem}

\newtheorem{lemma}[theorem]{Lemma}

\newtheorem{problem}[theorem]{Problem}
\newtheorem{proposition}[theorem]{Proposition}

\newenvironment{proof}[1][Proof]{\textbf{#1.} }{\ \rule{0.5em}{0.5em}}

\let\Section=\section
\def\section{\setcounter{equation}{0}\Section}

\begin{document}

\title{Stochastic integral representation of the  $L^{2}$ modulus of
 Brownian local time and a central limit theorem }
\author{Yaozhong Hu, David Nualart\thanks{
D. Nualart is supported by the NSF grant DMS0604207.} \\
Department of Mathematics \\
University of Kansas \\
Lawrence, Kansas, 66045 USA}
\date{}
\maketitle
\begin{abstract}
The purpose of this note is to prove a central limit theorem for the
$L^2$-modulus of continuity of the Brownian local time obtained in
\cite{CLMR}, using techniques of stochastic analysis. The main
ingredients of the proof are an asymptotic version of Knight's
theorem and the Clark-Ocone formula for the $L^2$-modulus  of the
Brownian local time.

\end{abstract}

\section{Introduction}

Let $B=\{B_{t},t\geq 0\}$ be a standard Brownian motion, and denote by $%
\{L_{t}(x),t\geq 0,x\in \mathbb{R\}}$ its local time. In \cite{CLMR} the
authors have proved the following central limit theorem for the $L^{2}$%
 modulus of continuity of the local time:

\begin{theorem}
\label{thm1}For each fixed $t>0$,%
\begin{equation}
h^{-\frac{3}{2}}\left( \int_{\mathbb{R}}(L_{t}(x+h)-L_{t}(x))^{2}dx-4th%
\right) \overset{\mathcal{L}}{\longrightarrow }8\sqrt{\frac{ \alpha_t}{3} }%
\eta ,  \label{b3}
\end{equation}%
as $h$ tends to zero, where
\begin{equation}
\alpha _{t}=\int_{\mathbb{R}}(L_{t}(x))^{2}dx,  \label{b5}
\end{equation}%
and $\eta $ is a $N(0,1)$ random variable independent of $B$.
\end{theorem}

We make use of the notation%
\begin{equation}
G_{t}(h)=\int_{\mathbb{R}}(L_{t}(x+h)-L_{t}(x))^{2}dx.  \label{b1}
\end{equation}%
It is proved in \cite[Lemma 8.1]{CLMR} that $\mathbb{E}\left(
G_{t}(h)\right) =4th+O(h^{2})$. Therefore, we can replace the term $4th$ in (%
\ref{b3}) by $\mathbb{\ E}\left( G_{t}(h)\right) .$

The proof of Theorem \ref{thm1}  is done in \cite{CLMR} by the
method of moments. The purpose of this paper is to provide a simple
proof of this result. Our method is based on an  asymptotic version
of Knight's theorem (see Revuz and Yor \cite{RY}, Theorem (2.3),
page 524) combined with the techniques of stochastic analysis and
Malliavin calculus.  The main idea is to apply the Clark-Ocone
stochastic integral representation \ formula to express $
\displaystyle G_{t}(h)\ -\mathbb{E}\left( G_{t}(h)\right) $ as a
stochastic integral. Then, by means of simple estimates using
H\"{o}lder's inequality, it is proved that the leading term is a
martingale,  to which we can apply an asymptotic version of Knight'
theorem. An important ingredient  is to show the convergence of the
quadratic variation of this martingale, which will be derived  by
using   Tanaka's formula  and backward It\^{o} stochastic integrals.

The paper is organized as follows. In the next section we recall some
preliminaries on Malliavin calculus and we establish a stochastic integral
representation for the random variable $G_{t}(h)$. Then, Section 3 is
devoted to the proof of Theorem \ref{thm1}.

\section{Stochastic integral representation of the $L^{2}$-modulus of
continuity}

Let us introduce some basic facts on the Malliavin calculus with respect the
the Brownian motion $B=\{B_{t},t\geq 0\}$. We refer \ to \cite{Nu06} for a
complete presentation of these notions. \ We assume that $B$ is defined on a
complete probability space $(\Omega ,\mathcal{F},P)$ such that $\mathcal{F}$
is generated by $B$. Consider the set $\mathcal{S}$ of smooth random
variables  of the form%
\begin{equation}
F=f\left( B_{t_{1}},\ldots ,B_{t_{n}}\right) ,  \label{h1}
\end{equation}%
where $\ t_{1},\ldots ,t_{n}\geq 0$, $f\in \mathcal{C}_{b}^{\infty }\left(
\mathbb{R}^{n}\right) $ (the space of bounded functions which have bounded
derivatives of all orders) and $n\in \mathbb{N}$. The derivative operator $D$
on a smooth random variable of the form (\ref{h1}) is defined by
\begin{equation*}
D_{t}F=\sum_{i=1}^{n}\frac{\partial f}{\partial x_{i}}\left(
B_{t_{1}},\ldots ,B_{t_{n}}\right) I_{[0,t_{i}]}(t),
\end{equation*}%
which is an element of $L^{2}\left( \Omega \times \lbrack 0,\infty )\right) $%
. We denote by \ $\mathbb{D}^{1,2}$ the completion of $\mathcal{S}$ with
respect to the norm $\left\| F\right\| _{1,2}\ $\ given by
\begin{equation*}
\left\| F\right\| _{1,2}^{2}=\mathbb{E}\left[ F^{2}\right] +\ \mathbb{E}%
\left( \int_{0}^{\infty }\left( D_{t}F\right) ^{2}dt\right) .
\end{equation*}%
The classical It\^{o} representation theorem asserts that any square
integrable random variable can be expressed as%
\begin{equation*}
F=\mathbb{E[}F\mathbb{]+}\int_{0}^{\infty }u_{t}dB_{t},
\end{equation*}
where $u=\{u_{t},t\geq 0\}$ is a unique adapted process such that $\mathbb{E}%
\left( \int_{0}^{\infty }u_{t}^{2}dt\right) <\infty $. If $F$ belongs to $%
\mathbb{D}^{1,2}$, then $u_{t}=\mathbb{E[D}_{t}F|\mathcal{F}_{t}]$, \ where $%
\{\mathcal{F}_{t},t\geq 0\}$ is the filtration generated by $B$, and we
obtain the Clark-Ocone formula (see \cite{Oc})%
\begin{equation}
F=\mathbb{E[}F\mathbb{]+}\int_{0}^{\infty }\mathbb{E[D}_{t}F|\mathcal{F}%
_{t}]dB_{t}.  \label{x2}
\end{equation}

The random variable
$G_{t}(h)=\int_{\mathbb{R}}(L_{t}(x+h)-L_{t}(x))^{2}dx$ can be
expressed in terms of the self-intersection local time of Brownian
motion. In fact,
\begin{eqnarray*}
&&G_{t}(h)=\int_{\mathbb{R}}\left[ \int_{0}^{t}\delta
(B_{u}+x+h)du-\int_{0}^{t}\delta (B_{u}+x)du\right] ^{2}dx \\
&=&2\left[ \int_{0}^{t}\int_{0}^{t}\delta
(B_{v}-B_{u}+h)dudv-\int_{0}^{t}\int_{0}^{t}\delta (B_{v}-B_{u})dudv\right]
\\
&=&2\int_{0}^{t}\int_{0}^{v}\left( \delta (B_{v}-B_{u}+h)+\delta
(B_{v}-B_{u}-h)-2\delta (B_{v}-B_{u})\right) dudv\,.
\end{eqnarray*}%
The rigorous justification of above argument can be made easily by
approximating the Dirac delta function by the heat kernel $p_{\varepsilon
}(x)=\frac{1}{\sqrt{2\pi \varepsilon }}e^{.x^{2}/2\varepsilon }$ as $%
\varepsilon $ tends to zero. That is, $G_{t}(h)$ is the limit in $%
L^{2}(\Omega )$ as $\varepsilon $ tends to zero of
\begin{equation}
G_{t}^{\varepsilon }(h)=2\int_{0}^{t}\int_{0}^{v}\left( p_{\varepsilon
}(B_{v}-B_{u}+h)+p_{\varepsilon }(B_{v}-B_{u}-h)-2p_{\varepsilon
}(B_{v}-B_{u})\right) dudv.  \label{x7}
\end{equation}%
Applying Clark-Ocone formula we can derive the following stochastic integral
representation for $G_{t}(h)$.

\begin{proposition}
\label{prop1}The random variable $G_{t}(h)$ defined in (\ref{b1}) can be
expressed as
\begin{equation*}
G_{t}(h)=\mathbb{E}(G_{t}(h))+\int_{0}^{t}u_{t,h}(r)dB_{r}\,,
\end{equation*}%
where
\begin{eqnarray}
u_{t,h}(r) &=&4\int_{0}^{r}\int_{0}^{h}\left( p_{t-r}(B_{r}-B_{u}+\eta
)-p_{t-r}(B_{r}-B_{u}-\eta )\right) d\eta du  \notag \\
&&+4\int_{0}^{r}\left( I_{[0,h]}(B_{r}-B_{u})-I_{[0,h]}(B_{u}-B_{r})\right)
du.  \label{e1}
\end{eqnarray}
\end{proposition}

\begin{proof}
For any $u<v$ and any $h\in \mathbb{R}$ we can write
\begin{equation*}
D_{r}p_{\varepsilon }(B_{v}-B_{u}+h)=p_{\varepsilon }^{\prime
}(B_{v}-B_{u}+h)I_{[u,v]}(r),
\end{equation*}%
and for any $u<r<v$%
\begin{eqnarray*}
\mathbb{E}\left( D_{r}p_{\varepsilon }^{\prime }(B_{v}-B_{u}+h)|\mathcal{F}%
_{r}\right)  &=&\mathbb{E}p_{\varepsilon }^{\prime }(\ \sqrt{v-r}\eta
+B_{r}-B_{u}+h) \\
&=&p_{v-r+\varepsilon }^{\prime }(B_{r}-B_{u}+h),
\end{eqnarray*}%
where $\eta $ denotes a $N(0,1)$ random variable independent of $B$.
Therefore, from Clark-Ocone formula (\ref{x2}) and Equation (\ref{x7}) we
obtain%
\begin{equation*}
G_{t}^{\varepsilon }(h)=\mathbb{E}(G_{t}^{\varepsilon
}(h))+\int_{0}^{t}u_{t,h}^{\varepsilon }(r)dB_{r},
\end{equation*}%
where%
\begin{eqnarray*}
u_{t,h}^{\varepsilon }(r) &=&2\int_{r}^{t}\int_{0}^{r}\Bigg(%
p_{v-r+\varepsilon }^{\prime }(B_{r}-B_{u}+h)+p_{v-r+\varepsilon }^{\prime
}(B_{r}-B_{u}-h) \\
&&-2p_{v-r+\varepsilon }^{\prime }(B_{r}-B_{u})\Bigg)dudv.
\end{eqnarray*}%
This expression can be written as
\begin{equation*}
u_{t,h}^{\varepsilon }(r)=2\int_{r}^{t}\int_{0}^{r}\int_{0}^{h}\left(
(p_{v-r+\varepsilon }^{\prime \prime }(B_{r}-B_{u}+\eta )-p_{v-r+\varepsilon
}^{\prime \prime }(B_{r}-B_{u}-\eta )\right) d\eta dudv.
\end{equation*}%
Using the fact that $p_{t}^{\prime \prime }(x)=2\frac{\partial p_{t}}{%
\partial t}(x)$ we obtain
\begin{eqnarray*}
u_{t,h}^{\varepsilon }(r) &=&4\left( \int_{0}^{r}\int_{0}^{h}\left(
(p_{t-r+\varepsilon }(B_{r}-B_{u}+\eta )-p_{t-r+\varepsilon
}(B_{r}-B_{u}-\eta )\right) d\eta du\right.  \\
&&\left. -\int_{0}^{r}\int_{0}^{h}\left( (p_{\varepsilon }(B_{r}-B_{u}+\eta
)-p_{\varepsilon }(B_{r}-B_{u}-\eta )\right) d\eta du\right) .
\end{eqnarray*}%
Letting $\varepsilon $ tend to zero we get \ that $u_{t,h}^{\varepsilon }(r)$
converges in $L^{2}(\Omega \times \lbrack 0,t])$ to $u_{t,h}(r)$ as $h$
tends to zero, which implies the desired result.
\end{proof}

\medskip From Proposition \ref{prop1} we can make the following
decomposition
\begin{equation*}
u_{t,h}(r)=\hat{u}_{t,h}(r)+\tilde{u}_{t,h}(r),
\end{equation*}%
where
\begin{eqnarray}
\hat{u}_{t,h}(r) &=&4\int_{0}^{r}\int_{0}^{h}\left( p_{t-r}(B_{r}-B_{u}+\eta
)-p_{t-r}(B_{r}-B_{u}-\eta )\right) d\eta du  \notag \\
&=&4\int_{0}^{r}\int_{0}^{h}\int_{-\eta }^{\eta }p_{t-r}^{\prime
}(B_{r}-B_{u}+\xi )d\xi d\eta du  \label{b4}
\end{eqnarray}%
and
\begin{equation*}
\tilde{u}_{t,h}(r)=4\int_{0}^{r}\left(
I_{[0,h]}(B_{r}-B_{u})-I_{[0,h]}(B_{u}-B_{r})\right) du.
\end{equation*}%
As a consequence, we finally obtain%
\begin{eqnarray}
G_{t}(h)-\mathbb{E}(G_{t}(h)) &=&4\int_{0}^{t}\hat{u}_{t,h}(r)dB_{r}  \notag
\\
&&\hspace{-0.24in}\hspace{-0.87in}+4\int_{0}^{t}\left( \int_{0}^{r}\left(
I_{[0,h]}(B_{r}-B_{u})-I_{[0,h]}(B_{u}-B_{r})\right) du\right) dB_{r}.
\label{x5}
\end{eqnarray}

\section{Proof of Theorem \ref{thm1}}

The proof will be done in several steps. Along the proof we will denote by $%
C $ a generic constant, which may be different form line to line.

\textit{Step 1 }We claim that the stochastic integral $\int_{0}^{t}$ $\hat{u}%
_{t,h}(r)dB_{r}$ makes no contribution to the limit  (\ref{b3}). That is,
\begin{equation*}
h^{-3/2}\int_{0}^{t}\hat{u}_{t,h}(r)dB_{r}
\end{equation*}%
converges in $L^{2}(\Omega )$ to zero as $h$ tends to zero. This is a
consequence of the next proposition.

\begin{proposition}
\label{prop2}There is a constant $C>0$ such that%
\begin{equation*}
\mathbb{E}\left( \int_{0}^{t}|\hat{u}_{t,h}(r)|^{2}dr\right) \leq Ch^{4}\,,
\end{equation*}%
for all $h>0$.
\end{proposition}

\begin{proof}
From (\ref{b4}) we can write%
\begin{eqnarray*}
\mathbb{E}\left( |\hat{u}_{t,h}(r)|^{2}\right)
&=&\int_{0}^{r}\int_{0}^{r}\int_{0}^{h}\int_{0}^{h}\int_{-\eta _{1}}^{\eta
_{1}}\int_{-\eta _{2}}^{\eta _{2}}\mathbb{E(}p_{t-r}^{\prime
}(B_{r}-B_{u_{1}}+\xi _{1}) \\
&&\times p_{t-r}^{\prime }(B_{r}-B_{u_{2}}+\xi _{2}))d\xi _{1}d\xi _{2}d\eta
_{1}d\eta _{2}du_{1}du_{2}.
\end{eqnarray*}%
By a symmetry argument, it suffices to integrate in the region $%
0<u_{1}<u_{2}<r$. Set
\begin{equation*}
\Phi (u_{1},u_{2},\xi _{1},\xi _{2})=\mathbb{E}\left( p_{t-r}^{\prime
}(B_{r}-B_{u_{1}}+\xi _{1})p_{t-r}^{\prime }(B_{r}-B_{u_{2}}+\xi
_{2})\right) .
\end{equation*}%
Then,
\begin{eqnarray*}
\Phi (u_{1},u_{2},\xi _{1},\xi _{2}) &=&\mathbb{E}\left( p_{t-r}^{\prime
}(B_{r}-B_{u_{2}}+B_{u_{2}}-B_{u_{1}}+\xi _{1})p_{t-r}^{\prime
}(B_{r}-B_{u_{2}}+\xi _{2})\right) \\
&=&\mathbb{E}\left( p_{t-r+u_{2}-u_{1}}^{\prime }(B_{r}-B_{u_{2}}+\xi
_{1})p_{t-r}^{\prime }(B_{r}-B_{u_{2}}+\xi _{2})\right) \\
&=&\int_{\mathbb{R}}p_{r-u_{2}}(z)p_{t-r+u_{2}-u_{1}}^{\prime }(z+\xi
_{1})p_{t-r}^{\prime }(z+\xi _{2})dz \\
&\leq &\Vert p_{r-u_{2}}\Vert _{p_{1}}\Vert p_{t-r+u_{2}-u_{1}}^{\prime
}\Vert _{p_{2}}\Vert p_{t-r}^{\prime }\Vert _{p_{3}}\,,
\end{eqnarray*}%
where $\frac{1}{p_{1}}+\frac{1}{p_{2}}+\frac{1}{p_{3}}=1$. It is easy to see
that
\begin{eqnarray*}
\Vert p_{r-u_{2}}\Vert _{p_{1}} &\leq &C(r-u_{2})^{-\frac{1}{2}+\frac{1}{%
2p_{1}}}, \\
\Vert p_{t-r+u_{2}-u_{1}}^{\prime }\Vert _{p_{2}} &\leq
&C(t-r+u_{2}-u_{1})^{-1+\frac{1}{2p_{2}}}\leq C(u_{2}-u_{1})^{-1+\frac{1}{%
2p_{2}}}, \\
\Vert p_{t-r}^{\prime }\Vert _{p_{3}} &\leq &C(t-r)^{-1+\frac{1}{2p_{3}}}\,,
\end{eqnarray*}%
for some constant $C>0$. Thus
\begin{eqnarray*}
\mathbb{E}\left( |\hat{u}_{t,h}(r)|^{2}\right) &\leq &C\
\int_{0}^{r}\int_{0}^{u_{2}}\int_{0}^{h}\int_{0}^{h}\int_{-\eta _{1}}^{\eta
_{1}}\int_{-\eta _{2}}^{\eta _{2}}(r-u_{2})^{-\frac{1}{2}+\frac{1}{2p_{1}}}
\\
&&\times (u_{2}-u_{1})^{-1+\frac{1}{2p_{2}}}(t-r)^{-1+\frac{1}{2p_{3}}}d\xi
_{1}d\xi _{2}d\eta _{1}d\eta _{2}du_{1}du_{2} \\
&\leq &C\ h^{4}\,.
\end{eqnarray*}%
This proves the proposition.
\end{proof}

\textit{Step 2 }Taking into account Proposition \ref{prop2}  and Equation (%
\ref{x5}), in order to show Theorem \ref{thm1} it suffices to show the
following convergence in law:
\begin{equation*}
h^{-\frac{3}{2}}\int_{0}^{t}\left( \int_{0}^{r}\left(
I_{[0,h]}(B_{r}-B_{u})-I_{[0,h]}(B_{u}-B_{r})\right) du\right) dB_{r}\overset%
{\mathcal{L}}{\rightarrow }2\eta \sqrt{\frac{\alpha _{t}}{3}\ },
\end{equation*}%
where $\eta $ is a standard normal random variable independent of $B$, and $%
\alpha _{t}$ has been defined in (\ref{b5}). Notice that
\begin{equation*}
M_{t}^{h}=h^{-\frac{3}{2}}\int_{0}^{t}\left( \int_{0}^{r}\left(
I_{[0,h]}(B_{r}-B_{u})-I_{[0,h]}(B_{u}-B_{r})\right) du\right) dB_{r}
\end{equation*}%
is a martingale with quadratic variation
\begin{equation*}
\left\langle M^{h}\right\rangle _{t}=h^{-3}\int_{0}^{t}\left(
\int_{0}^{r}\left( I_{[0,h]}(B_{r}-B_{u})-I_{[0,h]}(B_{u}-B_{r})\right)
du\right) ^{2}dr.
\end{equation*}%
From the asymptotic version of Knight's theorem (see Revuz and Yor \cite{RY}%
, Theorem (2.3) page. 524) it suffices to show the following convergences in
probability.
\begin{equation}
h^{-3}\int_{0}^{t}\left( \int_{0}^{r}\left(
I_{[0,h]}(B_{r}-B_{u})-I_{[0,h]}(B_{u}-B_{r})\right) du\right)
^{2}dr\rightarrow \frac{4}{3}\alpha _{t},  \label{a1}
\end{equation}%
and
\begin{equation}
\left\langle M^{h},B\right\rangle
_{t}=h^{-3/2}\int_{0}^{t}\int_{0}^{r}\left(
I_{[0,h]}(B_{r}-B_{u})-I_{[0,h]}(B_{u}-B_{r})\right) dudr\rightarrow 0,
\label{a2}
\end{equation}%
as $h$ tends to zero.
In fact, let $B^h$ be the Brownian motion such that  $ M^h _t = B^h_{\langle M^h \rangle_t}$.  Then, from
Theorem (2.3) pag. 524  in \cite{RY}, and the convergences  (\ref{a1}) and (\ref{a2}),
we deduce that $(B, B^h,  \langle M^h \rangle_t)$ converges in distribution to
$(B, \beta,  \frac 43 \alpha_t)$, where $\beta$ is a Brownian motion independent
of $B$. This implies that  $M_t^h=B^h_{\langle M^h \rangle_t}$ converges in
distribution to $\beta_ {\frac 43 \alpha_t}$, which yields the desired result.

Convergences  (\ref{a1}) and (\ref{a2}) will be proved in the following two
steps.

\textit{Step 3 }The convergence (\ref{a2}) follows from the following lemma.

\begin{lemma}
For any $t\geq 0$, $\left\langle M^{h},B\right\rangle _{t}$ converges to
zero in $L^{2}(\Omega )$ as $h$ tends to zero.
\end{lemma}

\begin{proof}
Fix $\varepsilon >0$. The quadratic covariation $\left\langle
M^{h},B\right\rangle _{t}$ is the limit, almost sure and in $L^{2}$, of
\begin{equation*}
\Phi _{\varepsilon
,h}=h^{-3/2}\int_{0}^{t}\int_{0}^{r}\int_{0}^{h}\int_{-\eta }^{\eta
}p_{\varepsilon }^{\prime }(B_{r}-B_{u}+\xi )d\xi d\eta dr
\end{equation*}%
as $\varepsilon $ tends to zero. We have%
\begin{eqnarray*}
\mathbb{E}\left( \Phi _{\varepsilon ,h}^{2}\right)
&=&h^{-3}\int_{0}^{t}\int_{0}^{t}\int_{0}^{r_{1}}\int_{0}^{r_{2}}%
\int_{0}^{h}\int_{0}^{h}\int_{-\eta _{1}}^{\eta _{1}}\int_{-\eta _{2}}^{\eta
_{2}}\mathbb{E}\big(p_{\varepsilon }^{\prime }(B_{r_{1}}-B_{u_{1}}+\xi _{1})
\\
&&p_{\varepsilon }^{\prime }(B_{r_{2}}-B_{u_{2}}+\xi _{2})\big)d\xi _{1}d\xi
_{2}d\eta _{1}d\eta _{2}du_{1}du_{2}dr_{1}dr_{2}.
\end{eqnarray*}%
We may assume that $0<r_{1}<r_{2}<t$, and we divide the region of
integration into three disjoint parts:

\textbf{Case 1}: $0<u_{1}<r_{1}<u_{2}<r_{2}<t$. In this situation
\begin{equation*}
\mathbb{E}\left[ p_{\varepsilon }^{\prime }(B_{r_{1}}-B_{u_{1}}+\xi
_{1})p_{\varepsilon }^{\prime }(B_{r_{2}}-B_{u_{2}}+\xi _{2})\right]
=p_{r_{1}-u_{1}+\varepsilon }^{\prime }(\xi _{1})p_{r_{2}-u_{2}+\varepsilon
}^{\prime }(\xi _{2}),
\end{equation*}%
which is an odd function of $\xi _{1}$ and $\xi _{2}$. Integration with
respect to $\xi _{1}$ from $-\eta _{1}$ to $\eta _{1}$ (or with respect to $%
\xi _{2}$ from $-\eta _{2}$ to $\eta _{2}$) yields $0$.

\textbf{Case 2}: $0<u_{1}<u_{2}<r_{1}<r_{2}<t$. In this case,
\begin{eqnarray*}
&&\mathbb{E}\left[ p_{\varepsilon }^{\prime }(B_{r_{1}}-B_{u_{1}}+\xi
_{1})p_{\varepsilon }^{\prime }(B_{r_{2}}-B_{u_{2}}+\xi _{2})\right]  \\
&=&\mathbb{E}\left[ p_{\varepsilon }^{\prime
}(B_{r_{1}}-B_{u_{2}}+B_{u_{2}}-B_{u_{1}}+\xi _{1})p_{\varepsilon }^{\prime
}(B_{r_{2}}-B_{r_{1}}+B_{r_{1}}-B_{u_{2}}+\xi _{2})\right]  \\
&=&\mathbb{E}\left[ p_{u_{2}-u_{1}+\varepsilon }^{\prime
}(B_{r_{1}}-B_{u_{2}}+\xi _{1})p_{r_{2}-r_{1}+\varepsilon }^{\prime
}(B_{r_{1}}-B_{u_{2}}+\xi _{2})\right]  \\
&=&\int_{\mathbb{R}}p_{r_{1}-u_{2}}(z)p_{u_{2}-u_{1}+\varepsilon }^{\prime
}(z+\xi _{1})p_{r_{2}-r_{1}+\varepsilon }^{\prime }(z+\xi _{2})dz\,.
\end{eqnarray*}%
The technique used in the proof of Proposition \ref{prop2} can be applied to
show that in this situation, we have%
\begin{eqnarray*}
&&h^{-3}\int_{0<u_{1}<u_{2}<r_{1}<r_{2}<t}\int_{0}^{h}\int_{0}^{h}\int_{-%
\eta _{1}}^{\eta _{1}}\int_{-\eta _{2}}^{\eta _{2}}\int_{\mathbb{R}%
}p_{r_{1}-u_{2}}(z) \\
&&\times \left| p_{u_{2}-u_{1}+\varepsilon }^{\prime }(z+\xi
_{1})p_{r_{2}-r_{1}+\varepsilon }^{\prime }(z+\xi _{2})\right| dzd\xi
_{1}d\xi _{2}d\eta _{1}d\eta _{2}du_{1}du_{2}dr_{1}dr_{2} \\
&\leq &Ct^{2}h\,.
\end{eqnarray*}

\textbf{Case 3}: $0<u_{2}<u_{1}<r_{1}<r_{2}<t$. In this case,
\begin{eqnarray*}
&&\mathbb{E}\left[ p_{\varepsilon }^{\prime }(B_{r_{1}}-B_{u_{1}}+\xi
_{1})p_{\varepsilon }^{\prime }(B_{r_{2}}-B_{u_{2}}+\xi _{2})\right]  \\
&=&\mathbb{E}\left[ p_{\varepsilon }^{\prime }(B_{r_{1}}-B_{u_{1}}+\xi
_{1})p_{\varepsilon }^{\prime
}(B_{r_{2}}-B_{r_{1}}+B_{r_{1}}-B_{u_{1}}+B_{u_{1}}-B_{u_{2}}+\xi _{2})%
\right]  \\
&=&\mathbb{E}\left[ p_{\varepsilon }^{\prime }(B_{r_{1}}-B_{u_{1}}+\xi
_{1})p_{r_{2}-r_{1}+u_{1}-u_{2}+\varepsilon }^{\prime
}(B_{r_{1}}-B_{u_{1}}+\xi _{2})\right]  \\
&=&\int_{\mathbb{R}}p_{r_{1}-u_{2}}(z)p_{\varepsilon }^{\prime }(z+\xi
_{1})p_{r_{2}-r_{1}+u_{1}-u_{2}+\varepsilon }^{\prime }(z+\xi _{2})dz.
\end{eqnarray*}%
Hence, in this case we have
\begin{eqnarray*}
&&h^{-3}\
\int_{0<u_{2}<u_{1}<r_{1}<r_{2}<t}\int_{0}^{h}\int_{0}^{h}\int_{-\eta
_{1}}^{\eta _{1}}\int_{-\eta _{2}}^{\eta _{2}}\int_{\mathbb{R}%
}p_{r_{1}-u_{2}}(z) \\
&&\times p_{\varepsilon }^{\prime }(z+\xi
_{1})p_{r_{2}-r_{1}+u_{1}-u_{2}+\varepsilon }^{\prime }(z+\xi _{2})dzd\xi
_{1}d\xi _{2}d\eta _{1}d\eta _{2}du_{1}du_{2}dr_{1}dr_{2} \\
&=&h^{-3}\ \int_{0<u_{2}<u_{1}<r_{1}<r_{2}<t}\int_{0}^{h}\int_{0}^{h}\int_{%
\mathbb{R}}p_{r_{1}-u_{2}}(z)\left[ p_{\varepsilon }(z+\eta
_{1})-p_{\varepsilon }(z-\eta _{1})\right]  \\
&&\times \left[ p_{r_{2}-r_{1}+u_{1}-u_{2}+\varepsilon }(z+\eta
_{2})-p_{r_{2}-r_{1}+u_{1}-u_{2}+\varepsilon }(z-\eta _{2})\right] dzd\eta
_{1}d\eta _{2}du_{1}du_{2}dr_{1}dr_{2}.
\end{eqnarray*}%
As $\varepsilon $ tends to zero, this converges to%
\begin{eqnarray*}
A(h) &:=&  2h^{-3}\
\int_{0<u_{2}<u_{1}<r_{1}<r_{2}<t}\int_{0}^{h}\int_{0}^{h}\
p_{r_{1}-u_{2}}(\eta _{1}) \\
&&\times \left[ p_{r_{2}-r_{1}+u_{1}-u_{2}}(\eta _{1}-\eta
_{2})-p_{r_{2}-r_{1}+u_{1}-u_{2}}(\eta _{1}+\eta _{2})\right] d\eta
_{1}d\eta _{2}du_{1}du_{2}dr_{1}dr_{2} \\
&=&-2h^{-3}\ \int_{0<u_{2}<u_{1}<r_{1}<r_{2}<t}\int_{0}^{h}\int_{0}^{h}\
\int_{-\eta _{2}}^{\eta _{2}}p_{r_{1}-u_{2}}(\eta _{1}) \\
&&\times p_{r_{2}-r_{1}+u_{1}-u_{2}}^{\prime }(\eta _{1}+\xi )d\xi d\eta
_{1}d\eta _{2}du_{1}du_{2}dr_{1}dr_{2}.
\end{eqnarray*}%
Notice that%
\begin{eqnarray*}
&&2h^{-3}\int_{0}^{h}\int_{0}^{h}\ \int_{-\eta _{2}}^{\eta
_{2}}p_{r_{1}-u_{2}}(\eta _{1})\left| p_{r_{2}-r_{1}+u_{1}-u_{2}}^{\prime
}(\eta _{1}+\xi )\right| d\xi d\eta _{1}d\eta _{2} \\
&\leq &Ch^{-3}\int_{0}^{h}\int_{0}^{h}\ \int_{-\eta _{2}}^{\eta
_{2}}(r_{1}-u_{2})^{-\frac{1}{2}}(r_{2}-r_{1}+u_{1}-u_{2})^{-1}d\xi d\eta
_{1}d\eta _{2} \\
&\leq &C(r_{1}-u_{2})^{-\frac{1}{2}}(r_{2}-r_{1}+u_{1}-u_{2})^{-1},
\end{eqnarray*}%
which is integrable in the region $\left\{
0<u_{2}<u_{1}<r_{1}<r_{2}<t\right\} $. \ On the other hand,%
\begin{eqnarray*}
&&\lim_{h\downarrow 0}2h^{-3}\int_{0}^{h}\int_{0}^{h}\ \int_{-\eta
_{2}}^{\eta _{2}}p_{r_{1}-u_{2}}(\eta
_{1})p_{r_{2}-r_{1}+u_{1}-u_{2}}^{\prime }(\eta _{1}+\xi )d\xi d\eta
_{1}d\eta _{2} \\
&=&2p_{r_{1}-u_{2}}(0)p_{r_{2}-r_{1}+u_{1}-u_{2}}^{\prime }(0)=0,
\end{eqnarray*}%
and, by the Dominated Convergence Theorem, we deduce that $A(h)$ converges
to zero as $h$ tends to zero. This completes the proof of the lemma.
\end{proof}

\textit{Step 4 }It only remains to show the convergence \ (\ref{a1}) in
probability. Define%
\begin{equation*}
\Psi _{h}(r)=\int_{0}^{r}\left(
I_{[0,h]}(B_{r}-B_{u})-I_{[0,h]}(B_{u}-B_{r})\right) du.
\end{equation*}%
We are going to show that%
\begin{equation}
h^{-3}\int_{0}^{t}\Psi _{h}(r)^{2}dr\overset{L^{2}(\Omega )}{\rightarrow }%
\frac{4}{3}\alpha _{t},  \label{f1}
\end{equation}%
as $h$ tends to zero.  Notice that
\[
\alpha_t= 2 \int_0^t\int_0^v  \delta_0(B_v -B_u) dudv
\]
is the self-intersection local time of $B$, and Equation (\ref{f1}) provides
an approximation for this self-intersection local time which has its own interest.

By the occupation formula for the Brownian motion we
can write%
\begin{eqnarray*}
\Psi _{h}(r) &=&\int_{\mathbb{R}}\left(
I_{[0,h]}(B_{r}-x)-I_{[0,h]}(x-B_{r})\right) L_{r}(x)dx \\
&=&\int_{0}^{h}\left( L_{r}(B_{r}-y)-L_{r}(B_{r}+y)\right) dy.
\end{eqnarray*}%
We can express the difference $L_{r}(B_{r}-y)-L_{r}(B_{r}+y)$ by means of
Tanaka's formula \ for the Brownian motion $\{B_{r}-B_{s},0\leq s\leq r\}$:%
\begin{eqnarray*}
L_{r}(B_{r}-y)-L_{r}(B_{r}+y) &=&-y-(B_{r}-y)^{+}+(B_{r}+y)^{+} \\
&&+\int_{0}^{r}\left( I_{B_{r}-B_{s}+y>0}-I_{B_{r}-B_{s}-y>0}\right) d%
\widehat{B}_{s},
\end{eqnarray*}%
where $d\widehat{B}_{s}$ denote the backward stochastic It\^{o} \ integral
and $y>0$. Integrating in the variable $y$ yields
\begin{eqnarray}
\Psi _{h}(r) &=&-\frac{h^{2}}{2}+\ \int_{0}^{h}\left[
(B_{r}+y)^{+}-(B_{r}-y)^{+}\right] dy  \notag \\
&&+\ \int_{0}^{h}\left( \int_{0}^{r}I_{y>|B_{r}-B_{s}|}d\widehat{B}%
_{s}\right) dy.  \label{b8}
\end{eqnarray}%
\ By stochastic Fubini's theorem%
\begin{equation}
\int_{0}^{h}\left( \int_{0}^{r}I_{y>|B_{r}-B_{s}|}d\widehat{B}_{s}\right)
dy=\int_{0}^{r}(h-|B_{r}-B_{s}|)^{+}d\widehat{B}_{s}.\   \label{b9}
\end{equation}%
Taking into account  (\ref{b8}) and (\ref{b9}),  the convergence (\ref{f1})
\ will follow from%
\begin{equation}
\ h^{-3}\int_{0}^{t}\left( \int_{0}^{r}(h-|B_{r}-B_{s}|)^{+}d\widehat{B}%
_{s}\right) ^{2}dr\overset{L^{2}(\Omega )}{\rightarrow }\frac{4}{3}\alpha
_{t},  \label{f2}
\end{equation}%
as $h$ tends to zero. By It\^{o}'s formula we can write%
\begin{eqnarray}
&&\left( \int_{0}^{r}(h-|B_{r}-B_{s}|)^{+}d\widehat{B}_{s}\right) ^{2}  \notag
 = 2\int_{0}^{r}\left( \int_{s}^{r}(h-|B_{r}-B_{u}|)^{+}d\widehat{B}%
_{u}\right)  \\
&&\qquad \times  (h-|B_{r}-B_{s}|)^{+}d\widehat{B}_{s}  \label{w1}
  +\int_{0}^{r}\left[ (h-|B_{r}-B_{s}|)^{+}\right] ^{2}ds.
\end{eqnarray}%
Finally, (\ref{f2}) follows form (\ref{w1}) and the next two lemmas.

\begin{lemma}
\label{Lem1}We have%
\begin{equation*}
\ \int_{0}^{t}\int_{0}^{r}\frac{\left[ (h-|B_{r}-B_{s}|)^{+}\right] ^{2}}{%
h^{3}}dsdr\overset{L^{2}(\Omega )}{\rightarrow }\frac{4}{3}\alpha _{t},
\end{equation*}%
as $h$ tends to zero.
\end{lemma}

\begin{proof}
Notice that%
\begin{equation*}
\alpha _{t}=\int_{\mathbb{R}}L_{t}(x)^{2}dx=\int_{0}^{t}\int_{\mathbb{R}%
}L_{r}(x)L_{dr}(x)dx=\int_{0}^{t}L_{r}(B_{r})dr,
\end{equation*}%
and%
\begin{equation*}
\ \int_{0}^{t}\int_{0}^{r}\frac{\left[ (h-|B_{r}-B_{s}|)^{+}\right] ^{2}}{%
h^{3}}dsdr=\ \int_{0}^{t}\int_{\mathbb{R}}\frac{\left[ (h-|B_{r}-x|)^{+}%
\right] ^{2}}{h^{3}}L_{r}(x)dx.
\end{equation*}%
As a consequence, taking into account that
\[
\int_{\mathbb{R}}\frac{\left[
(h-|B_{r}-x|)^{+}\right] ^{2}}{h^{3}}dx=  \int_{\mathbb{R}} \frac {[(h-|x|)^+]^2} {h^3} dx =
\frac{4}{3},
\]
 we obtain
\begin{eqnarray*}
&&\left| \int_{0}^{t}\int_{0}^{r}\frac{\left[ (h-|B_{r}-B_{s}|)^{+}\right]
^{2}}{h^{3}}dsdr-\frac{4}{3}\alpha _{t}\right| \\
&\leq &\int_{0}^{t}\int_{\mathbb{R}}\frac{\left[ (h-|B_{r}-x|)^{+}\right]
^{2}}{h^{3}}\left| L_{r}(x)-L_{r}(B_{r})\right| dxdr \\
&\leq &\frac{4}{3}\int_{0}^{t}\sup_{|x-y|<h}\left| L_{r}(x)-L_{r}(y)\right|
dr,
\end{eqnarray*}%
which clearly converges to zero in $L^{2}$ by the properties of the Brownian
local time (see, for instance, \cite{B}).
\end{proof}

\begin{lemma}
We have%
\begin{equation*}
\frac{1}{h^{6}}\mathbb{E}\left[ \left( \int_{0}^{t}\left( \int_{0}^{r}\left(
\int_{s}^{r}(h-|B_{r}-B_{u}|)^{+}d\widehat{B}_{u}\right)
(h-|B_{r}-B_{s}|)^{+}d\widehat{B}_{s}\right) dr\right) ^{2}\right]
\rightarrow 0,
\end{equation*}%
as $h$ tends to zero.
\end{lemma}

\begin{proof}
By the isometry property of the backward It\^{o} integral we can write%
\begin{eqnarray*}
B_{h} &:&=\frac{1}{h^{6}}\mathbb{E}\left[ \left( \int_{0}^{t}\left(
\int_{0}^{r}\left( \int_{s}^{r}(h-|B_{r}-B_{u}|)^{+}d\widehat{B}_{u}\right)
(h-|B_{r}-B_{s}|)^{+}d\widehat{B}_{s}\right) dr\right) ^{2}\right]  \\
&=&\frac{1}{h^{6}}\mathbb{E}\left[ \int_{0}^{t}\left(
\int_{s}^{t}(h-|B_{r}-B_{s}|)^{+}\left( \int_{s}^{r}\ (h-|B_{r}-B_{u}|)^{+}d%
\widehat{B}_{u}\right) dr\right) ^{2}ds\right]  \\
&\leq &\ \mathbb{E}\left[ B_{h}^{1}B_{h}^{2}\ \right] ,
\end{eqnarray*}%
where%
\begin{equation*}
B_{h}^{1}=\int_{0}^{t}\left( \int_{s}^{t}\frac{(h-|B_{r}-B_{s}|)^{+}}{h^{2}}%
\ dr\right) ^{2}ds
\end{equation*}%
and%
\begin{equation*}
B_{h}^{2}=\sup_{0<s<r<t}\left| \int_{s}^{r}\ \frac{(h-|B_{r}-B_{u}|)^{+}}{h}d%
\widehat{B}_{u}\right| ^{2}.
\end{equation*}%
As in Lemma \ref{Lem1}, we can show that $B_{h}^{1}$ converges to \ $\frac{9%
}{4}\int_{0}^{t}\left( L_{t}(B_{s})-L_{s}(B_{s})\right) ^{2}ds$, and the
convergence holds in $L^{p}(\Omega )$ for any $p\geq 2$.  Here we use
\[
\int_{\mathbb{R}}  \frac {(h-|x|)^+} {h^2} dx = \frac 32.
\]
 On the other
hand, $\int_{s}^{r}\ \frac{(h-|B_{r}-B_{u}|)^{+}}{h}d\widehat{B}_{u}$ can be
expressed using again Tanaka's formula:%
\begin{eqnarray*}
\frac{1}{h}\int_{s}^{r}\ (h-|B_{r}-B_{u}|)^{+}d\widehat{B}_{u}
&=&\int_{0}^{h}\left( \int_{s}^{r}I_{y>|B_{r}-B_{u}|}d\widehat{B}_{u}\right)
dy \\
&=&\frac{1}{h}\int_{0}^{h}\left( L_{r}(B_{r}-y)-L_{r}(B_{r}+y)\right) dy \\
&&-\frac{1}{h}\int_{0}^{h}\left( L_{s}(B_{r}-y)-L_{s}(B_{r}+y)\right) dy+%
\frac{h}{2} \\
&&+\ \frac{1}{h}\int_{0}^{h}\left[ (B_{r}-B_{s}+y)^{+}-(B_{r}-B_{s}-y)^{+}%
\right] dy.
\end{eqnarray*}%
Therefore,$\ $%
\begin{equation*}
\left| \frac{1}{h}\int_{s}^{r}\ (h-|B_{r}-B_{u}|)^{+}d\widehat{B}_{u}\right|
\leq \sup_{s<r<t}\sup_{|x-y|\leq 2h}\ \left| L_{r}(x)-L_{r}(y)\right| +O(h),
\end{equation*}%
which also converges to zero in \ $L^{p}(\Omega )$ for any $p\geq 2$.
\end{proof}

\end{document}